\documentclass[12pt, reqno]{elsarticle}

\usepackage{amsmath, amssymb, amsfonts, amsbsy, amsthm, latexsym}
\usepackage{graphicx}
\usepackage{verbatim}
\usepackage{color}
\usepackage{algorithm}
\usepackage{algorithmic}
\usepackage{caption}
\usepackage{subcaption}

\newtheorem{theorem}{Theorem}[section]

\newtheorem{lemma}[theorem]{Lemma}

\theoremstyle{definition}
\newtheorem{definition}[theorem]{Definition}

\numberwithin{equation}{section}

\numberwithin{equation}{section}

\usepackage{enumerate}

\def\0{{\bar{0}}}

\DeclareMathOperator*{\argmin}{\arg\!\min}

\begin{document}

\title{An Alternative Thresholding Rule for Compressed Sensing}

\author{Jonathan Ashbrock}
\address{Department of Mathematics, Vanderbilt University, Nashville, TN 37240, USA}
\begin{abstract}
	Compressed Sensing algorithms often make use of the hard thresholding operator to pass from dense vectors to their best $s-$sparse approximations. However, the output of the hard thresholding operator does not depend on any information from a particular problem instance. We propose an alternative thresholding rule, Look Ahead Thresholding, that does. In this paper we offer both theoretical and experimental justification for the use of this new thresholding rule throughout compressed sensing.
\end{abstract}

\maketitle



\

\section{Introduction}

Real world signals often contain relatively small amounts of information compared to their dimension. This property manifests in the following widely observed fact: many signal classes, when expanded in the appropriate basis, are very sparse. Sparse signals inherently are simpler to store, transmit, and process than dense ones. Because of the utility of sparse signals, it is of particular interest to find sparse representations in useful bases when they exist.

The field of Compressed Sensing (CS) was developed precisely to study the concepts of representation, acquisition, and recovery of sparse signals. One central question in CS is this: Given a vector of possibly noisy measurements, how may we find the sparsest signal that explains these measurements (up to, perhaps, some error tolerance level). More precisely, letting  $\|x\|_0$ denote the number of non-zero elements of the vector $x$, then given a tolerance $\epsilon$, we wish to solve \begin{align}
\min_{x \in \mathbb{C}^d} \|x\|_0 \hspace{5mm} \text{ subject to } \|Ax-y\|_2 <\epsilon \label{CS_Problem}
\end{align}
when given both $A$ and $y$.

While the field of CS is mainly concerned with finding the sparse representations, it is necessary to study closely related questions. For instance, it is well known that \eqref{CS_Problem} is NP-Hard for general matrices $A$ and every $\epsilon \geq 0$ \cite{Foucart}. However, with additional structure forced upon $A$, namely the restricted isometry property (RIP),  the problem is actually tractable. 

With this knowledge that the problem \textit{can} be solved, research quickly turned to the question of \textit{how} to solve \eqref{CS_Problem}. The list of techniques includes, yet is not limited to, iterative hard thresholding (IHT) \cite{IHT_for_Sparse,IHT}, (orthogonal) matching pursuit \cite{OMP}, hard thresholding pursuit \cite{HTP}, and CoSaMP \cite{COSAMP}. Each of these algorithms solve \eqref{CS_Problem} provided $A$ has sufficiently strong RIP.

Many of the above techniques, particularly those which are \textit{thresholding} based, rely crucially on the action of the hard thresholding operator. For instance, IHT simply alternates a gradient descent step followed by hard thresholding. The hard thresholding operator acts by taking in an arbitrary vector and returning the nearest $s-$sparse vector. However, doing so makes no use of $y$, the measurement vector. This leads us to ask: Can we find a better thresholding technique?

In this paper, we answer this question in the affirmative. We propose the Look Ahead Thresholding (LAT) technique which will be shown in a variety of settings to achieve better results than simply using hard thresholding. In particular, we modify IHT to use our new thresholding rule and propose the Iterative Look Ahead Thresholding (ILAT) algorithm.  We prove that ILAT has comparable worst case performance to IHT in terms of the required RIP to converge. Moreover, we show both experimentally and theoretically that look ahead thresholding excels when used in compressed sensing.

The remainder of the paper is organized as follows. Section \ref{background_section} will present the background information, develop our new thresholding rule, and present the ILAT algorithm. Section \ref{noiseless_section} analyzes the worst-case behavior of ILAT in relation to the RIP of the sensing matrix. Section \ref{average_section} contains an average case analysis showing the power of look ahead thresholding. Finally, Section \ref{experiments_section} shows in a few different experiments that this technique performs exceedingly well in practice.

\section{Background Information}\label{background_section}
Let us begin with the necessary priors from compressed sensing while simultaneously standardizing our notation. A vector $x\in \mathbb{R}^d$ is said to be $s-$sparse if it has at most $s$ non-zero entries. We note that most signals of interest are only sparse in a particular basis so one may similarly define sparsity with respect to some other orthonormal basis such as wavelets, fourier, etc. 

The process of `taking a measurement' is computing an inner product against a known vector. So, the transformation from signal to a vector of measurements is realized by matrix multiplication. We let $x^*\in \mathbb{R}^d$ be the $s-$sparse signal to recover. Given access to the vector of measurements $y=Ax^*$ in  $\mathbb{R}^m$, the goal is to recover quickly and uniquely the $s-$sparse vector $x^*$. Notice that if $m\geq d$ and $A$ is invertible, we can recover $x^*$ perfectly from $y$. In practice, we want to minimize the necessary number of measurements so we study at length the case when $m<d$. 

At the other extreme, notice that if we knew the support of $x^*$, then $s$ linearly independent measurements would suffice to recover $x^*$. While we cannot know the support of every measured signal, knowing that the signal is sparse inherently reduces the complexity of $x^*$. Therefore, we hope to still be able to get away with using $m$ measurements with $s<m<d$. 

Two of the main questions, then, in compressed sensing are: 
\begin{enumerate}
	\item What does the matrix $A$ need to look like for $x^*$ to be uniquely determined amongst all $s-$sparse vectors by $y=Ax^*$?
	\item Knowing $A$ and $y$, how can we procedurally reconstruct the solution $x^*$?
\end{enumerate}

The first question is very well-studied and traditionally answered by enforcing $A$ to have the restricted isometry property. 

\begin{definition}
	A matrix $A$ has $s^{th}$ restricted isometry constant $\delta_s=\delta_s(A)$ provided that $\delta_s$ is the smallest $\delta \geq 0$ so that for every $s-$sparse vector $x$ we have\begin{align*}
	(1-\delta)\|x\|_2^2 \leq \|Ax\|_2^2 \leq (1+\delta)\|x\|_2^2.
	\end{align*}
\end{definition}
We say informally that $A$ has the restricted isometry property (RIP) if $\delta_s$ is relatively small for $s$ relatively large. Essentially, though $A$ cannot possibly be an isometry when $m<d$, if $A$ has RIP of order $s$ then every subset of $s$ columns of $A$ form an approximate isometry. 

We let $\|\cdot \|_{op}$ be the operator norm of a matrix relative to the Euclidean norms on both the domain and range. Every matrix $A$ trivially has $\delta_s(A)\leq \max\{1,\|A\|_{op}^2-1\}$ so we need to ask:  How strong does the RIP of $A$ need to be to guarantee that $x^*$ is uniquely determined by $y$? In theory, recovery is possible if the matrix $A$ is injective on the set of $s-$sparse vectors. Notice that if $x,y$ are both $s-$sparse, their difference is $2s-$sparse. Then, if $\delta_{2s}<1$, it follows that $A(x-y)\neq 0$ by definition. So, in fact the condition $\delta_{2s}<1$ is sufficient for the map $x^*\mapsto Ax^*$ to have \textit{an} inverse on the $s-$sparse signal space.

While $\delta_{2s}<1$ is sufficient \textit{in theory} for recovery of $s-$sparse vectors, in practice this is not enough. Notice that the set of sparse vectors is not a subspace of $\mathbb{R}^d$, so the inverse map $Ax^*\mapsto x^*$ is not linear. Thus, most compressed sensing techniques require much stronger versions of the restricted isometry property in order to be able to find $x^*$ quickly, robustly, and stably. 

Unfortunately, the restricted isometry property is a delicate thing. Constructing matrices of a given size and with suitable RIP is difficult or impossible to do deterministically. Luckily, though, large enough random matrices with entries chosen i.i.d. from a Bernoulli or other sub-Gaussian distribution tend to have sufficiently strong RIP \cite{Foucart}. That is, we do not know necessarily how to construct RIP matrices but we do know that using random matrices with sufficiently many rows will suffice with high probability. 

Having answered question $(1)$ above, we can now turn our attention to question $(2)$: how exactly can we perform the non-linear inversion of the matrix $A$ on the set of $s-$sparse vectors? There are countless techniques and algorithms, one of which we discuss in depth in the next section.

\subsection{Iterative Hard Thresholding}
Iterative Hard Thresholding was first introduced by Blumensath and Davies in \cite{IHT_for_Sparse} as a method of sparse dictionary approximation and by the same authors shortly after in \cite{IHT} for compressed sensing. There are many ways to interpret IHT but we prefer to understand it as alternating a gradient descent step followed by a hard thresholding step. 

To develop the gradient descent setting, given the measurements $y$ and matrix $A$, define the cost function \begin{align}
C(x)=\|y-Ax\|_2^2 = \|A(x^*-x)\|_2^2 \label{cost_function}
\end{align} to describe how well the vector $x$ explains the observed measurements $y$. A routine calculation shows the gradient of $C$ is $\nabla C(x)=-2A^*(y-Ax)=-2A^*A(x^*-x)$. In general, a thresholder is any function that outputs \textit{some} $s-$sparse projection of the input vector. In particular, we define $H_s$ to be the \textit{hard thresholding operator} which acts on a vector by retaining only the $s$ largest entries in magnitude while setting the rest to zero (ties can be broken arbitrarily). Then, given a starting point $x^0$ and using a `step-size' of $1/2$, IHT is defined by the iteration\begin{enumerate}
	\item $a^{t+1} = x^t-\frac{1}{2}\nabla C(x^t)=x^t+A^*(y-Ax^t)$
	\item $x^{t+1}=H_s(a^{t+1})$.
\end{enumerate}

The inspiration here is that gradient descent helps move towards a zero of $C$ while the thresholding ensures that our iterates are actually sparse. And, under the conditions $\delta_{2s}<1$ the only point having both these properties is $x^*$. However, because it is possible that $H_s$ `undoes' enough of the progress of the gradient step, we are not necessarily guaranteed to converge to a point where $\nabla C(x^t)=0$. If we specify a stronger RIP, though, we can guarantee convergence. 

For IHT to work, we need $x^t$ to stay \textit{about as close} to $x^*$ as $a^t$ is. That is, $H_s$ needs to not undo too much of the gradient progress. The following crucially important fact from \cite{IHT} gets us most of the way there:
 \begin{align}
\|x^*-x^{t}\|_2 =\|x^*-H_s(a^t)\|_2\leq 2 \|x^*-a^{t}\|_2\label{IHT_bound}.
\end{align}
Then, because $x^*-a^{t}=(I-A^*A)(x^*-x^{t-1})$ and noting that the RIP forces $A^*A$ to be close to the identity when applied to sparse vectors (see Lemma \ref{RIP_lemma}), we can show $\|x^*-x^{t}\|_2 \leq \rho \|x^*-x^{t-1}\|_2$ for $\rho<1$ so that $(x^t)$ converges to $x^*$. 

This fact generalizes. Suppose we have an algorithm that alternates $a^t=x^{t-1}+A^*(y-Ax^{t-1})$ then chooses $x^t$ to be \textit{some} thresholding of $a^t$ (not necessarily using $H_s$). Then, if \begin{align}
\|x^*-x^t\|_2\leq k\|x^*-a^t\|_2\label{k_inequality}
\end{align}
for some universal $k$, the same Lemma \ref{RIP_lemma} lets us specify which RIP of $A$ (namely $\delta_{2s}<1/k^2$) suffices to prove convergence when normalizing $\|A\|_{op}=1$. This is the proof technique we pursue in Section \ref{noiseless_section}.

Here is a final thought to motivate our campaign for alternative thresholding. If $P_*$ is the projection onto the support of $x^*$, then \begin{align*}
\|x^*-P_* a^t\|_2\leq \|x^*-a^t\|_2.
\end{align*}
That is, inequality \eqref{k_inequality} is satisfied for $k=1$. Comparing this to the hard thresholding requirement $k=2$, this leaves quite a bit of room for improvement. While we cannot use $P_*$ because the support is in general unknown, maybe we can use some of the information content from the measurements $y$ to pick a thresholding that does better than $k=2$. 
\subsection{Look Ahead Thresholding}
In this section we will describe how our novel look ahead thresholder chooses which coordinates to keep and which to kill. Along the way, we will show the thought process that led to the definition of look ahead thresholding.

Fix a vector $z\in \mathbb{R}^d$. Then, $H_s(z)$ is the $s-$sparse vector $H(z)$ minimizing $\|z-H(z)\|_2$. However, per the discussion preceding inequality \eqref{k_inequality}, the ideal thresholder $H$ for sparse recovery is the one which minimizes $\|x^*-H(z)\|_2$. That is, it is preferable to threshold a vector $z$ in a way so that the result is as close to the solution $x^*$ as possible. While we cannot minimize $\|x^*-H(z)\|_2$ directly because we do not know $x^*$, we develop look ahead thresholding to get a good approximation. 

Notice that the value of the cost $C$ at a particular thresholding $H(z)$ of $z$ is \begin{align*}
C(H(z))=\|y-A(H(z))\|_2^2=\|A(x^*-H(z))\|_2^2.
\end{align*}Then, because $x^*-H(z)$ is $2s-$sparse, if $A$ has RIP it is an approximate isometry so that $C(H(z))$ is close to $\|x^*-H(z)\|_2^2$. Moreover, $C$ is something we can work with unlike $\|x^*-H(z)\|_2^2$.

 Unfortunately, $C$ is a non-separable quadratic and there are $\binom{d}{s}$ possible $s-$sparse projections of a point $z$. It is computationally intractable to find the $H(z)$ minimizing $C$. We address this issue by introducing a surrogate $f_{\eta,z}$ which is a separable quadratic related closely to the Taylor series of $C$ at $z$. The idea for using $f_{\eta,z}$ was inspired by the work in \cite{Optimal_Brain_Damage} which searched for sparse representations for neural networks.

First let us simply define the surrogate \begin{align}
f_{\eta, z}(x)=C(z)+ \sum_{i=1}^d \left(2\eta \frac{\partial C}{\partial x_i}(z)(x_i-z_i) + (x_i-z_i)^2\right)\label{weighted_surrogate_def}.
\end{align}
We claim this function is quite closely related to the Taylor series for $C$ at $z$. Indeed, first assume we have normalized the columns of $A$ to have $\|A_i\|_2=1$ so that $\frac{\partial^2 C}{\partial x_i^2}=1$. Then the function $f_{\eta,z}$ is formed from the Taylor series for $C$ by first ignoring the off-diagonal terms of the Hessian then multiplying the linear terms by a tunable constant. Ignoring the off-diagonal terms is necessary for separability and adding the weight to the gradient term lets us control the algorithm's performance more carefully.

Now we can define the look ahead thresholding rule. If $z$ is the point to threshold, we  pick the thresholding of $z$ that minimizes $f_{\eta, z}$. More precisely, first let $P_\Omega$ be the projection onto the coordinate axes indexed by $\Omega\subset \{1,\dots, d\}$.  Then, the look ahead thresholding function $H_{s,\eta}$ has the action defined by
\begin{align}
H_{s,\eta}(z)=\argmin_{P_\Omega z, |\Omega|=s}f_{\eta,z}(P_\Omega z). \label{diagonal_hessian_projection}
\end{align}

We recall that the original problem \eqref{CS_Problem} can be solved exactly if given unlimited time. Therefore, it is important to emphasize that the action defined in \eqref{diagonal_hessian_projection} can be computed quickly. Notice that, for a particular projection $P_\Omega(z)$, the function $f_{\eta, z}$ takes the value \begin{align}
f_{\eta, z}(P_\Omega(z))=C(z)+\sum_{i\notin \Omega} \left(2\eta \frac{\partial C}{\partial x_i}(z)(-z_i) +  (-z_i)^2\right).\label{projection_surrogate_value}
\end{align}

Then we observe that the solution to \eqref{diagonal_hessian_projection} can be computed by picking $\Omega$ to contain precisely the $s$ coordinates for which the values \begin{align}
2\eta\frac{\partial C}{\partial x_i}(z)(-z_i)+(z_i)^2\label{LAT_Rule}
\end{align}
are largest. Because computation of $\frac{\partial C}{\partial x_i}(z)$ is by far the most expensive step of \eqref{LAT_Rule}, the cost to compute $H_{s,\eta} (z)$ is essentially the cost of computing the gradient $\nabla C(z)$.

Now, we promised a second motivation for the definition of the thresholding rule \eqref{diagonal_hessian_projection}. This reformulation of $H_{s,\eta}$ is the basis for the name \textit{look ahead thresholding}. Because this reformulation is used throughout the proofs in Sections \ref{noiseless_section} and \ref{average_section}, Lemma \ref{reformulation_lemma} proves this reformulation is correct. Moreover, this reformulation also lets us see that hard thresholding is a special case of look ahead thresholding with the weight $\eta$ chosen to be zero. 

 Define the \textbf{look ahead point} $\ell_\eta=z-\eta\nabla C(z)$, the point that would be the next gradient descent step from $z$ with step-size $\eta$. Then, look ahead thresholding defined by \eqref{diagonal_hessian_projection} picks the $s-$sparse projection of $z$ that is closest to $\ell_\eta$.
Because for specific values of $\eta$, the point $\ell_\eta$ is provably closer to $x^*$ than $z$ is, we might expect this thresholding to find better projections. Indeed, this is the case.

To show how look ahead thresholding may be applied in practice, we define the Iterative Look Ahead Thresholding algorithm as an example. The only difference between ILAT and IHT is that the thresholding step is done using $H_{s,\eta}$ instead of $H_s$. That is, we use look ahead thresholding instead of hard thresholding to hopefully retain as much of the progress made by the gradient updates as possible while remaining $s-$sparse.

In Section \ref{noiseless_section} we show that this algorithm is guaranteed to converge under some RIP conditions that are, perhaps surprisingly, more restrictive than the conditions for IHT. However, convergence guarantees are worst case analyses and we argue for the use of look ahead thresholding both theoretically in Section \ref{average_section} and experimentally in Section \ref{experiments_section} by appealing to the average case.  Let us end this section by formally presenting Iterative Look Ahead Thresholding.

\begin{algorithm}
	\textbf{Iterative Look Ahead Thresholding}
\begin{algorithmic}
		\STATE\textbf{Input: } Matrix $A$, measurements $y$, sparsity $s$,  iterations $T$
\STATE \textbf{Output: } Estimate $x^T$ of $s-$sparse solution $x^*$.	
\STATE Set $x^0=0\in \mathbb{R}^d$
	 \FOR{$1 \leq t \leq T$}
	 \STATE Set $a^{t}=x^{t-1}+A^*(y-Ax^{t-1})$
	 \FOR{$1\leq i \leq d$}
	 \STATE Compute $s_i=2\eta(-a^{t}_i) \frac{\partial C}{\partial x_i}(a^t) +(a^t_i)^2$
	 \ENDFOR
	 \STATE Set $\Omega \subset \{1,\dots, d\}$ the indices of the $s$ largest values of $s_i$. 
	 \STATE Set $x^t=P_\Omega a^t$.
	 \ENDFOR
	\end{algorithmic}
\end{algorithm}

\section{Iterative Look Ahead Thresholding Converges with Suitable RIP}\label{noiseless_section}

As is typical for compressed sensing, we present theorems of the form `for sufficient restricted isometry constants, our algorithm is guaranteed to recover the solution $x^*$'. We prove theorems of this sort both in a noiseless environment and a corresponding result when the measurements may be corrupted. We begin first with a few technical facts that show up in our analysis time and again.

\subsection{Technical Lemmas}
 We point out that in our analysis, we require that the measurement matrix $A$ is normalized to have operator norm $1$. This is slightly different from traditional CS literature which normalizes the columns of $A$ to have expected squared length $1$. However, in some settings this normalization has been observed to yield superior performance \cite{Normalized_IHT}. 
\begin{lemma}\label{eigenvalue_lemma}
	If $A$ is an $m\times d$ matrix and $\|A\|_{op}=1$, then $\|I-2\eta A^*A\|_{op}\leq 1$ whenever $\eta \in [0,1]$. Moreover, if $m<d$ then $\|I-2\eta A^*A\|_{op}=1$.
\end{lemma}
\begin{proof}
	We build a singular value decomposition for $I-2\eta A^*A$ from an svd for $A$ itself. Let the svd for $A$ be $A=V\Sigma U$. Then notice that $I-2\eta A^*A = U^*(I-2\eta \Sigma^*\Sigma)U$. Notice first that because $\|A\|_{op}=1$, the diagonal of $\Sigma^*\Sigma$ is in $[0,1]$. Thus, the diagonal of $D=I-2\eta \Sigma^*\Sigma$ is in $[-1,1]$. Now, form $\tilde{D}$ and $\tilde{U}$ by first replacing the negative entries in $D$ by their absolute values then multiplying the corresponding row in $U$ by $-1$. 
	
	We point out three things. First, that we can write $I-2\eta A^*A = U^*\tilde{D}\tilde{U}$. Second, the matrix $\tilde{U}$ remains a unitary. Finally, $\tilde{D}$ is a diagonal matrix with positive entries in $[0,1]$. Therefore $U\tilde{D}\tilde{U}$ is an svd for $I-2\eta A^*A$ so $\|I-2\eta A^*A\|_{op}\leq 1$. 
	
	 Now, because $m<d$ we can pick a nontrivial $v$ in the kernel of $A$ so $(I-2\eta A^*A)v=v$. Therefore $\|I-2\eta A^*A\|_{op}=1$. 
\end{proof}

The prior lemma controlled the operator norm of $I-2\eta A^*A$ over the entirety of its domain. The following lemma controls the operator norm of $I-A^*A$ when restricted to $s-$sparse vectors. In this case we get much better performance in terms of the restricted isometry constants of $A$.
\begin{lemma}\label{RIP_lemma}
	If $x\in \mathbb{R}^d$ is $s-$sparse and $\|A\|_{op}=1$, then $\|(I-A^*A)x\|_2 \leq \sqrt{\delta_s}\|x\|_2$.
\end{lemma}

\begin{proof}
	Simply expand \begin{align}
	\|(I-A^*A)x\|_2^2&=\|x\|_2^2 + \|A^*Ax\|_2^2-2\langle x, A^*Ax\rangle \nonumber\\
	&= \|x\|_2^2 -\|Ax\|_2^2 + \|A^*Ax\|_2^2 - \|Ax\|_2^2 \label{first_RIP_ineq}.
	\end{align}
	Now, because $A^*$ has norm $1$,  $\|A^*(Ax)\|_2^2 -\|Ax\|_2^2\leq 0$. By re-arranging the RIP condition we may also see \begin{align}
	\|Ax\|_2^2 \geq (1-\delta_s)\|x\|_2^2 \Rightarrow \delta_s\|x\|_2^2\geq \|x\|_2^2-\|Ax\|_2^2 \label{second_RIP_ineq}.
	\end{align}
	Then, substituting the reformulation in \eqref{second_RIP_ineq} into \eqref{first_RIP_ineq}, our result follows by taking a square root. 
\end{proof}

We remark that Lemma \ref{RIP_lemma} is similar to the following reformulation of $\delta_s$ whose proof can be found in \cite{Foucart}. If $A_\Omega$ is the matrix formed from $A$ by taking the columns indexed by $\Omega$, then \begin{align}
\delta_s(A) = \sup_{|\Omega|=s}\|I-A_\Omega^*A_\Omega\|_{op}^2. \label{RIP_reformulation}
\end{align}
Finally, we will show that using the decision rule \eqref{diagonal_hessian_projection} to threshold a vector $z$ is equivalent to picking the $s-$sparse thresholding of $z$ which is closest to the look ahead point $\ell_\eta$. 
\begin{lemma}\label{reformulation_lemma}
	Fix a vector $z\in \mathbb{R}^d$, $\eta>0$, measurement vector $y\in \mathbb{R}^m$, and an $m\times d$ matrix $A$. Let $C(z)=\|y-Az\|_2^2$ and define $\ell_\eta = z-\eta\nabla C(z)$. Then, the vector $H_{s,\eta} (z)$ is the closest $s-$sparse thresholding of $z$ to $\ell_\eta$. 
\end{lemma}

\begin{proof}
	Notice that for a generic vector $x$, \begin{align}
	\|x-\ell_\eta\|_2^2 &= \sum_{i=1}^d \left( x_i -\left(z_i-\eta\frac{\partial C}{\partial x_i}(z)\right)\right)^2 \nonumber\\
	&=\sum_{i=1}^d \left((x_i-z_i)+\eta\frac{\partial C}{\partial x_i}(z)\right)^2\nonumber\\
	&=\sum_{i=1}^d \left(2\eta\frac{\partial C}{\partial x_i}(z)(x_i-z_i) +  \left(x_i-z_i\right)^2\right) +\eta^2\|\nabla C(z)\|_2^2\label{look_ahead_taylor_equivalence}.
	\end{align}
	For an arbitrary index set $\Omega$ and projection $P_\Omega z$ we have that \begin{align*}
	\|P_\Omega z-\ell^t_\eta\|_2^2 = \sum_{i\notin \Omega}\left(2\eta\frac{\partial C}{\partial x_i}(z)(-z_i) +  \left(z_i\right)^2\right)+\eta^2\|\nabla C(z)\|_2^2.
	\end{align*}
	Now notice that, no matter the choice of $\Omega$, the value \begin{align*}
	\sum_{i\notin \Omega}\left(2\eta\frac{\partial C}{\partial x_i}(z)(-z_i) +  \left(z_i\right)^2\right)+\sum_{i\in \Omega}\left(2\eta\frac{\partial C}{\partial x_i}(z)(-z_i) +  \left(z_i\right)^2\right)+\eta^2\|\nabla C(z)\|_2^2
	\end{align*}
	is a constant. So, minimizing the sum indexed over $i\notin \Omega$ is equivalent to maximizing the sum indexed over $i\in \Omega$. Selecting $\Omega$ in this way is precisely the decision rule for $H_{s,\eta} (z)$ given by \eqref{LAT_Rule}.
\end{proof}

\subsection{Noiseless Convergence}

In this section we prove that Iterative Look Ahead Thresholding recovers exactly the exactly $s-$sparse solution $x^*$ from the noiseless measurements $y=Ax$ provided $\delta_{2s}(A)$ is sufficiently small. Per the discussion preceding \eqref{k_inequality}, if we can establish a universal $k$ so that  \begin{align}\|x^*-x^t\|_2 =\|x^*-H_{s,\eta}(a^t)\|_2\leq k\|x^*-a^t\|_2\label{goal_bound}\end{align} then invoking Lemma \ref{RIP_lemma} and expanding the definition of $a^t$, it follows that \begin{align*}
\|x^*-x^t\|_2 \leq k\sqrt{\delta_{2s}}\|x^*-x^{t-1}\|_2.
\end{align*}

 Then, we may specify an RIP constraint so that we recover $x^*$ in the limit. The first step is to determine a suitable constant $k$.

\begin{lemma}\label{Projection_lemma}
	Suppose we are given an $n\times d$ matrix $A$ with $\|A\|_{op}=1$ and $y=Ax^*$ for some $s-$sparse $x^*\in \mathbb{R}^d$. Let $0\leq \eta \leq 1$ and the sequences $(x^t),(a^t)$ be defined by Iterative Look Ahead Thresholding with parameter $\eta$. Then, \begin{align*}
	\|x^*-x^t\|_2\leq (1+\sqrt{1+4\eta^2})\|x^*-a^t\|_2.
	\end{align*}
\end{lemma}
\begin{proof}
	Let $\ell_\eta=a^t-\eta \nabla C(a^t)$, the look ahead point from $a^t$. Let $P$ be the projection onto the support of a best $s-$term approximation of $\ell_\eta$. Because $x^t=H_{s,\eta}(a^t)$ is the closest $s-$sparse thresholding of $a^t$ to $\ell_\eta$, it follows that $\|x^t-\ell_\eta\|_2\leq \|P a^t -\ell_\eta\|_2$. Then we write \begin{align}\|x^*-x^t\|_2&\leq\|x^*-\ell_\eta\|_2+ \|x^t-\ell_\eta\|_2\nonumber\\
	&\leq \|x^*-\ell_\eta\|_2+\|P a^t-\ell_\eta\|_2 \label{factorized_form}.
	\end{align}
	Notice first that $\|x^*-\ell_\eta\|_2=\|(I-2\eta A^*A)(x^*-a^t)\|_2$. Because $0\leq \eta \leq 1$ and through Lemma \ref{eigenvalue_lemma} we can say \begin{align}\|x^*-\ell_\eta \|_2\leq \|x^*-a^t\|_2\label{look_ahead_closer}.\end{align} 
	Dealing now with the second term of \eqref{factorized_form} we see
	\begin{align}
	\|P a^t -\ell_\eta\|_2^2 &= \|P a^t-P \ell_\eta\|_2^2 + \|P \ell_\eta - \ell_\eta\|_2^2\nonumber \\
	&\leq \|a^t-\ell_\eta\|_2^2+\|x^*-\ell_\eta\|_2^2\label{first_step}\\
	&\leq \|a^t-\ell_\eta\|_2^2+\|x^*-a^t\|_2^2
	\label{projection_bound}.
	\end{align} 
	Above, inequality \eqref{first_step} is because projection is norm $1$ and because $P \ell_\eta$ is a better $s-$term approximation of $\ell_\eta$ than $x^*$ is. Now again because of the normalization $\|A\|_{op}=1$, \begin{align*}\|a^t-\ell_\eta\|_2 = \|\eta \nabla C(a^t)\|_2=\|2\eta A^*A(x^*-a^t)\|_2\leq2\eta\|x^*-a^t\|_2.\end{align*} We finish working with inequality \eqref{projection_bound} by substituting to yield: \begin{align}
	\|P a^t-\ell_\eta\|_2^2 \leq (1+4\eta^2) \|x^*-a^t\|_2^2.\label{last_projection_bound}
	\end{align}
	Finally, substituting \eqref{look_ahead_closer} and \eqref{last_projection_bound} into \eqref{factorized_form} we see \begin{align*}
	\|x^t-x^*\|_2\leq \|x^*-a^t\|_2+\sqrt{1+4\eta^2}\|x^*-a^t\|_2
	\end{align*}
	which is precisely the inequality given in the Lemma statement.
\end{proof}

Now, as argued above, the proof of our first main theorem requires just an algebraic manipulation after appealing to Lemmas \ref{Projection_lemma} and \ref{RIP_lemma}.
\begin{theorem}\label{noiseless_theorem}
	Suppose we are given an $m\times d$ matrix $A$ with $\|A\|_{op}=1$ and $y=Ax^*$ for some $s-$sparse $x^*\in \mathbb{R}^d$. Let $0\leq \eta \leq 1$ and the sequences $(x^t),(a^t)$ be defined by Iterative Look Ahead Thresholding using with parameter $\eta$. If the restricted isometry constant $\delta_{2s}$ of $A$ satisfies $\delta_{2s}<\frac{1}{\left(1+\sqrt{1+4\eta^2}\right)^2},$ then \begin{align*}
	\|x^t-x^*\|_2 \leq \rho^t \|x^0-x^*\|_2
	\end{align*}
	where $\rho=\sqrt{\delta_{2s}}\left(1+\sqrt{1+4\eta^2}\right)<1$. In particular, $(x^t)$ converges to $x^*$.
\end{theorem}

\begin{proof}
		Because the conditions of Lemma \ref{Projection_lemma} are satisfied, we have that $\|x^t-x^*\|_2\leq (1+\sqrt{1+4\eta^2})\|a^t-x^*\|_2$. Then,\begin{align*}
	\|x^t-x^*\|_2&\leq (1+\sqrt{1+4\eta^2})\|a^t-x^*\|_2\\
	&=(1+\sqrt{1+4\eta^2})\|x^{t-1}-x^*+A^*A(x^*-x^{t-1})\|_2\\
	&=(1+\sqrt{1+4\eta^2})\|(I-A^*A)(x^{t-1}-x^*)\|_2\\
	&\leq (1+\sqrt{1+4\eta^2})\sqrt{\delta_{2s}}\|x^{t-1}-x^*\|_2
	\end{align*}
	where the final inequality follows from Lemma \ref{RIP_lemma} because $x^{t-1}-x^*$ is $2s$-sparse and $\|A\|_{op}=1$. Therefore, by iteration we have that $\|x^t-x^*\|_2\leq \rho^t \|x^0-x^*\|_2$ for $\rho<1$ so convergence is guaranteed.
\end{proof}

\subsection{Noisy Convergence}

Let us now turn our attention to the case where there is noise present. In particular, the signal $x^*$ is no longer required to be exactly $s-$sparse and the measurements look like $y=Ax^*+e$ for $e$ some noise term. We define the following auxiliary noise term which describes how much $y$ is corrupted from the noiseless measurements of the best $s-$term approximation of $x^*$.

\begin{definition}
	For a vector $x^* \in \mathbb{R}^d$, $H_s(x^*)$ its best $s-$term approximation, and a corrupted measurement vector $y=Ax^* + e$, we define the error term $\tilde{e} =y-A(H_s(x^*))=A(x^*-H_s(x^*))+e$.
\end{definition}

We have the following useful relation for the gradient at a point in terms of $\tilde{e}$: \begin{align}
\eta\nabla C(z) = -2\eta A^*(y-Az)
&=-2\eta A^*(Ax^*+e-Az)\nonumber\\
&=-2\eta A^*A(H_s(x^*)-z)+A^*(A(x^*-H_s(x^*))+e)\nonumber\\
&= -2\eta A^*A(H_s(x^*)-z)-2\eta A^*\tilde{e}\label{noise_rearrangement}.
\end{align}

Before proving our noisy convergence theorem, which proceeds very similarly to the previous section, let us offer a comment on what convergence actually means here. Because the signal $x^*$ is not exactly $s-$sparse and our algorithm returns sparse signals, at best we hope to recover the signal $H_s(x^*)$. However, because our measurements are corrupted by noise, even this cannot be expected. We do show, however, that $(x^t)$ will converge to some small neighborhood of $H_s(x^*)$ and the size of the neighborhood depends only on $\delta_{2s},\|\tilde{e}\|_2$ and $\eta$.

\begin{lemma}\label{noisy_projection_bound}
	Suppose we are given an $m\times d$ matrix $A$ with $\|A\|_{op}=1$ and $y=Ax^*+e$ for some $x^*\in \mathbb{R}^d$. Let $0\leq\eta \leq 1$ and the sequences $(x^t),(a^t)$ be defined by Iterative Look Ahead Thresholding with parameter $\eta$. Then, \begin{align*}
	\|x^t-H_s(x^*)\|_2\leq (2+2\eta)\|a^t-H_s(x^*)\|_2+6\eta \|\tilde{e}\|_2.
	\end{align*}
\end{lemma}

\begin{proof}
	Let $\ell_\eta=a^t-\eta \nabla C(a^t)$. First we will derive two bounds based on the equality \eqref{noise_rearrangement}. Because Lemma \ref{eigenvalue_lemma} guarantees $\|I-2\eta A^*A\|_{op}\leq 1$, 
	\begin{align}
	\|H_s(x^*)-\ell_\eta\|_2=\|H_s(x^*)-a^t+\eta \nabla C(a^t)\|_2&=\|(I-2\eta A^*A)(H_s(x^*)-a^t)+2\eta A^*\tilde{e}\|_2\nonumber\\
	&\leq \|H_s(x^*)-a^t\|_2 + 2\eta \|\tilde{e}\|_2 \label{look_ahead_length_part1}.
	\end{align}
	Second, directly from \eqref{noise_rearrangement} we have 
	\begin{align}
	\|a^t-\ell_\eta\|_2=\|\eta\nabla C(a^t)\|_2 \leq 2\eta \|H_s(x^*)-a^t\|_2+2\eta \|\tilde{e}\|_2\label{second_look_ahead_length}.
	\end{align}
Next we use the triangle inequality: \begin{align}\|x^t-H_s(x^*)\|_2\leq \|x^t-\ell_\eta\|_2+\|H_s(x^*)-\ell_\eta\|\label{triangle_inequality}.
	\end{align} Because the second term on the right of \eqref{triangle_inequality} is bounded by \eqref{look_ahead_length_part1}, let us work with the first term on the right. If $P$ is the projection onto the support of $H_s(\ell_\eta)$, then because $x^t$ is the closest $s-$sparse thresholding of $a^t$ to $\ell_\eta$,  we also have that \begin{align}
	\|x^t-\ell_\eta\|_2\leq \|Pa^t-\ell_\eta\|_2&\leq \|Pa^t-P\ell_\eta\|_2+\|P\ell_\eta -\ell_\eta\|_2\nonumber \\
	&\leq \|a^t-\ell_\eta\|_2 + \|H_s(x^*)-\ell_\eta\|_2,\label{random_label}
	\end{align}
	where the last term is because $P\ell_\eta$ is a better $s-$term approximation of $\ell_\eta$ than is $H_s(x^*)$. Substituting \eqref{random_label} into \eqref{triangle_inequality} yields \begin{align}
	\|x^t-H_s(x^*)\|_2 \leq 2\|H_s(x^*)-\ell_\eta\|_2+\|a^t-\ell_\eta\|_2\label{almost_done}.
	\end{align}
	
	Now, we again substitute \eqref{look_ahead_length_part1} and \eqref{second_look_ahead_length} directly into \eqref{almost_done} achieves the stated bound in this lemma.
\end{proof}

Again, our second main theorem follows from Lemma \ref{noisy_projection_bound} with some simple algebra. 

\begin{theorem}\label{noisy_theorem}
	Suppose we are given an $m\times d$ matrix $A$ with $\|A\|_{op}=1$ and $y=Ax^*+e$ for $x^*\in \mathbb{R}^d$. Let $0\leq\eta \leq 1$ and the sequences $(x^t),(a^t)$ be defined by Iterative Look Ahead Thresholding with parameter $\eta$. If the restricted isometry constant $\delta_{2s}$ of $A$ satisfies $\delta_{2s}<\frac{1}{(2+2\eta)^2},$ then the sequence $(x^t)$ satisfies the recurrence relation \begin{align*}
	\|H_s(x^*)-x^t\|_2\leq \rho \|H_s(x^*)-x^{t-1}\|_2+(2+8\eta)\|\tilde{e}\|_2
	\end{align*}
	for $\rho = (2+2\eta)\sqrt{\delta_{2s}}<1$.
	\end{theorem}
\begin{proof}
	By Lemma \ref{noisy_projection_bound}, we have $
	\|H_s(x^*)-x^t\|_2\leq (2+2\eta)\|H_s(x^*)-a^t\|_2+6\eta\|\tilde{e}\|_2$.
	But notice that \begin{align}
		\|H_s(x^*)-x^t\|_2&\leq (2+2\eta)\|H_s(x^*)-a^t\|_2+6\eta\|\tilde{e}\|_2\nonumber\\
		&=(2+2\eta)\|H_s(x^*)-x^{t-1}-A^*(y-Ax^{t-1})\|_2+6\eta\|\tilde{e}\|_2\nonumber\\
		&=(2+2\eta)\|H_s(x^*)-x^{t-1}-A^*(A(H_s(x^*)-x^{t-1}+x^*-H_s(x^*))+e)\|_2+6\eta\|\tilde{e}\|_2\nonumber\\
		&\leq (2+2\eta)\left(\|(I-A^*A)(H_s(x^*)-x^{t-1})\|_2+\|A^*(\tilde{e})\|_2\right)+6\eta\|\tilde{e}\|_2\nonumber\\
		&\leq (2+2\eta)\sqrt{\delta_{2s}}\|H_s(x^*)-x^{t-1}\|_2+(2+8\eta)\|\tilde{e}\|_2\label{noisy_RIP_bound}
		\end{align}
		where the last line is because of Lemma \ref{RIP_lemma}. 
\end{proof}
\section{Average Case Analysis}\label{average_section}

Our goal in defining a better thresholding rule $H$ given an $s-$sparse vector to recover $x^*$ is to minimize the distance $\|H(z)-x^*\|_2$. Moreover, we are interested in measuring the size of $\|H(z)-x^*\|_2$ relative to $\|z-x^*\|_2$. Recall that in the case of $H=H_s$, the hard thresholding operator, at best we can say that $\|H_s(z)-x^*\|_2\leq 2\|z-x^*\|_2$ (see, e.g., inequality $(21)$ of \cite{IHT}). Even if we take the expected value of $\|H_s(z)-x^*\|$ over a random choice of measurement matrix $A$, this upper bound does not improve because $H_s$ is independent of $A$.

However, for look ahead thresholding the average case is much better. In this section we will show that taking expectation over random $A$ with entries chosen i.i.d. from a Gaussian, \begin{align}
\mathbb{E}\left[\| H_{s,\eta} (z)-x^*\|_2\right] \leq \rho \|z-x^*\|_2
\end{align}
and that $\rho$ can be taken strictly smaller than $2$ for a range of $\eta$ values. Now, during the proof of Lemma \ref{Projection_lemma} we derive the following important inequality which is the basis of this section's analysis. We recreate it here in the precise form we need for completeness.

\begin{lemma}\label{transformation_lemma}
	Let $z\in \mathbb{R}^d$ be a vector, $x^*\in \mathbb{R}^d$ be $s-$sparse and $A$ some $m\times d$ measurement matrix. Let $C(z)=\|A(x^*-z)\|_2^2$ and define $\ell_\eta = z-\eta \nabla C(z)$ to be the look ahead point from $z$. Then, \begin{align}
	\|H_{s,\eta} (z)-x^*\|_2\leq \sqrt{\|z-\ell_\eta\|_2^2+\|\ell_\eta-x^*\|_2^2} +\|\ell_\eta-x^*\|_2. \label{analysis_inequality}
	\end{align}
\end{lemma}

\begin{proof}
	Let $P$ be the projection onto the support of a best $s-$term approximation of $\ell_\eta$. Then, because $x^*$ is a worse $s-$term approximation of $\ell_\eta$ than $P\ell_\eta$ is, \begin{align*}
	\|Pz-\ell_\eta\|_2^2 &= \|Pz-P\ell_\eta\|_2^2 + \|\ell_\eta-P\ell_\eta\|_2^2\\
	&\leq \|z-\ell_\eta\|_2^2 + \|\ell_\eta-x^*\|_2^2.
	\end{align*}
	Now, by choice $H_{s,\eta}(z)$ is at least as close to $\ell_\eta$ as $Pz$ is. So, a triangle inequality leaves \begin{align*}
	\|H_{s,\eta} (z)-x^*\|_2 &\leq \|H_{s,\eta}(z)-\ell_\eta\|_2 + \|\ell_\eta-x^*\|_2\\
	&\leq \|Pz-\ell_\eta\|_2+\|\ell_\eta-x^*\|_2\\
	&\leq \sqrt{\|z-\ell_\eta \|_2^2+\|\ell_\eta-x^*\|_2^2} +\|\ell_\eta - x^*\|_2.
	\end{align*}
\end{proof}

The reason this upper bound is useful is because we can write each term on the right hand side of \eqref{analysis_inequality} as a transformation of a common vector:\begin{align}
z-\ell_\eta &= -2\eta A^*A(x^*-z) \label{first_random_vec}\\
x^*-\ell_\eta &= (I-2\eta A^*A)(x^*-z) \label{second_random_vec}.
\end{align}
Then, controlling the random behavior of $A$ and the related matrices $I-2\eta A^*A$ and $2\eta A^*A$ gives us our result. In the rest of this section we compute the expected size of $\|z-\ell_\eta\|_2^2$ and $\|x^*-\ell_\eta\|_2^2$ over the random draw of $A$. 

We first desire to answer: For random $V$, in what ways does the average size of $\|Vy\|_2^2$ depend on both $\|y\|_2^2$ and on $V$? It turns out that under some mild constraints which are satisfied in our setting (see Lemma \ref{randomness_verification}), the important quantity is the expected Frobenius norm of the matrix in question. 

\begin{lemma}\label{Frobenius_Lemma}
	Fix a vector $y\in \mathbb{R}^d$ and let $V$ be an $m\times d$ random matrix with entries satisfying \begin{enumerate}
		\item The columns $V_i, 1\leq i \leq d$ of $V$ all have the same expected squared length.
		\item Two distinct entries $V_{i,j},V_{i,k}$ from the same row have $\mathbb{E}\left[V_{i,j}V_{i,k}\right]=0$.
	\end{enumerate}
Then, the expected value of $\|Vy\|_2^2$ is scaled by the Frobenius norm like:\begin{align*}
\mathbb{E}\left[ \|Vy\|_2^2\right] = \frac{1}{d}\cdot \mathbb{E}\left[\|V\|_F^2\right]\cdot \|y\|_2^2.
\end{align*}
\end{lemma}
\begin{proof}
	We begin by expanding the random variable of concern:  \begin{align}
	\|Vy\|_2^2 &=\sum_{i=1}^m \left( \sum_{j=1}^d V_{i,j}y_j\right) \left(\sum_{k=1}^d V_{i,k}y_k\right)\nonumber\\
	&= \sum_{i=1}^m\sum_{j=1}^d\sum_{k=1}^d V_{i,j}V_{i,k}y_jy_k\label{squared_norm_expansion}.
	\end{align}
	
	If $j\neq k$, then $\mathbb{E}\left[ V_{i,j}V_{i,k}y_jy_k\right]=0$ because of the assumption $(2)$ on $V$. Therefore taking the expectation of both sides of \eqref{squared_norm_expansion} we are only left with \begin{align}
	\mathbb{E}\left[ \|Vy\|_2^2\right] &= \sum_{i=1}^m\sum_{j=1}^d \mathbb{E}\left( V_{i,j}^2y_j^2\right) \nonumber\\
	&=\sum_{j=1}^d y_j^2\sum_{i=1}^m \mathbb{E}\left(V_{i,j}^2\right)\label{middle_squared_norm_expansion}.
	\end{align}
	
	Then, by the assumption $(1)$ on the random matrix $V$, the terms $\sum_{i=1}^m\mathbb{E}(V_{i,j}^2)$ are constant for each column $j$. Letting this common squared expected length be denoted $c$, then we can continue \eqref{middle_squared_norm_expansion} by \begin{align}
	\mathbb{E}\left[ \|Vy\|_2^2\right] = c\sum_{j=1}^dy_j^2=c\|y\|_2^2 \label{expectation}.
	\end{align}
	Finally, it remains to show that the the constant $c$ is equal to the expected squared Frobenius norm divided by the dimension of $y$. Well, the expected Frobenius norm can be written \begin{align*}
	\mathbb{E}\left[ \|V\|_F^2\right] =\sum_{j=1}^d\sum_{i=1}^m V_{i,j}^2=\sum_{j=1}^dc=cd
	\end{align*}
	which completes the proof
\end{proof}

We remark that the probabilistic setting above is quite mild and is satisfied by most random matrix ensembles including any matrix with entries selected i.i.d. from a zero mean distribution. For the remainder of this section we will work with $A$ whose entries are drawn i.i.d. from a Gaussian normalized so the expected squared column length is $1$. Before computing the expected Frobenius norms, though, we will isolate a few key probabilistic facts. Throughout the remainder of this paper we let $N(\mu, \sigma^2)$ denote the normal distribution with mean $\mu$ and variance $\sigma^2$. 

\begin{lemma}\label{inner_product_lemma}
	Let $A_1,A_2,A_3\in \mathbb{R}^m$ be random vectors with entries selected $i.i.d.$ from $N(0,\frac{1}{m})$. Then, \begin{align}
	\mathbb{E}\left[ \langle A_1,A_2\rangle^2\right]=\frac{1}{m} \label{inner_product_squared}
	\end{align}
	and for $p=3$ or $p=1$ we have \begin{align}
	\mathbb{E}\left[ \langle A_1, A_2\rangle\langle A_1, A_p\rangle \right] =0\label{different_inner_products}.
	\end{align}
	Moreover, the fourth moment of $\|A_i\|_2$ behaves like \begin{align}
		\mathbb{E}\left[ \|A_i\|_2^4\right] =1+\frac{2}{m}.
	\end{align}
\end{lemma}
\begin{proof}
	We expand the first expectation like \begin{align}
	\mathbb{E}\left( \langle A_1, A_2\rangle ^2\right) &= \mathbb{E}\left[\left(\sum_{k=1}^m (A_1)_k(A_2)_k\right)\cdot\left(\sum_{\ell=1}^m (A_1)_\ell (A_2)_\ell\right)\right]\nonumber\\
	&=\mathbb{E} \left[ \sum_{k,\ell=1}^m (A_1)_k (A_2)_k(A_1)_\ell (A_2)_\ell\right].\label{inner_product_expansion}
	\end{align}
	First notice that when $k\neq \ell$, the four terms are distinct and chosen mutually independently so expectation commutes with taking their product. Then, since they are each zero mean, the product of the four terms has zero expected value. Therefore, after taking the expectation of both sides of \eqref{inner_product_expansion} the only terms contributing are $k=\ell$. Therefore, \begin{align*}
	\mathbb{E}\left[ \langle A_1,A_2\rangle ^2\right] &=\sum_{k=1}^m \mathbb{E}\left[ (A_1)_k^2(A_2)_k^2\right]\\
	&= \sum_{k=1}^m \mathbb{E}\left[(A_1)_k^2\right] \mathbb{E}\left[(A_2)_k^2\right] = \frac{1}{m}
	\end{align*}
	where the last line is because $(A_1)_k$ and $(A_2)_k$ are independent. 
	
	Now for the second part of this lemma. Identically to the previous part, \begin{align*}
	\mathbb{E}\left[ \langle A_1,A_2\rangle\langle A_1,A_p\rangle\right] = \sum_{k,\ell=1}^m \mathbb{E}\left[(A_1)_k(A_2)_k(A_1)_\ell (A_p)_\ell\right].
	\end{align*}
	Now, because the entries are mutually independent and $(A_2)_k$ is always distinct from $(A_1)_k,(A_1)_\ell$ and $(A_p)_\ell$ no matter if $p=1$ or $p=3$, then for every $k,\ell$ we have \begin{align*}
	\mathbb{E}\left[ (A_1)_k(A_2)_k(A_1)_\ell (A_p)_\ell\right]=\mathbb{E}\left[(A_1)_k (A_1)_\ell (A_p)_\ell\right]\mathbb{E}\left[ (A_2)_k\right]=0
	\end{align*}
	which finishes the second part of this lemma.
	
	Finally we will compute the fourth moment of $\|A_j\|_2$ \begin{align}
	\|A_j\|_2^4 = \left(\sum_{i=1}^m (A_j)_i^2\right) \left(\sum_{k=1}^m (A_j)_i^2\right)\nonumber \\
	=\sum_{i=1}^m (A_j)_k^4 + \sum_{i=1}^m\sum_{k\neq i}(A_j)_i^2(A_j)_k^2 \label{fourth_moment_expanded}.
	\end{align}
	Now because $(A_j)_i,(A_j)_k$ are independent when $i\neq k$, taking expectations of both sides of \eqref{fourth_moment_expanded} leaves us with \begin{align}
	\mathbb{E}\left[ \|A_j\|_2^4\right] &=\sum_{i=1}^m \mathbb{E}\left[ (A_j)_i^4\right] + \sum_{i=1}^m\sum_{k\neq j}\mathbb{E}\left[ (A_j)_i^2\right] \mathbb{E}\left[(A_j)_k^2\right]\nonumber\\
	&=m\left(\frac{3}{m^2}\right)+(m^2-m)\left(\frac{1}{m^2}\right)\nonumber
	\end{align}
	where the last equality is in part because the fourth moment of a zero mean Gaussian is $3\sigma^4$. Cleaning up this expression finishes the proof. 
\end{proof}

With all our tools ready, we can now compute the expected Frobenius norm of the matrix $I-2\eta A^*A$. 

\begin{lemma}\label{I_mins}
	Let $A$ be an $m\times d$ random matrix with entries $A_{i,j}$ selected i.i.d. from $N(0,\frac{1}{m})$ and let $\eta >0$. Then the squared Frobenius norm of $I-2\eta A^*A$ has expected value \begin{align*}
	\mathbb{E}\left[ \|I-2\eta A^*A\|_F^2\right] = 4d\left[ \left( \frac{d+m+1}{m}\right)\eta^2 -\eta +\frac{1}{4}\right].
	\end{align*}
\end{lemma}

\begin{proof}
	The off diagonal entries of $I-2\eta A^*A$, of which there are $d^2-d$, all have the form $-2\eta \langle A_i,A_j\rangle$ for $i\neq j$. Then, by Lemma \ref{inner_product_lemma} the expected squared size of each off diagonal term is $\frac{4\eta^2}{m}$. Therefore the off-diagonal contribution to the squared Frobenius norm is $(d^2-d)\frac{4\eta^2}{m}$. 
	
	The diagonal entries of $I-2\eta A^*A$ are of the form $1-2\eta \|A_i\|^2_2$. So their squared expectation again by the Lemma \ref{inner_product_lemma} is\begin{align*}
	\mathbb{E}\left[ (1-2\eta\|A_i\|_2^2)^2\right] &= 1-4\eta \mathbb{E}\left[\|A_i\|_2^2\right] + 4\eta^2\mathbb{E}\left[ \|A_i\|_2^4\right] \\
	&=1-4\eta + 4\eta^2\left(1+\frac{2}{m}\right)
	\end{align*}
	so the diagonal contribution to the expected squared Frobenius norm is $d\left( 1-4\eta + 4\eta^2\left(1+\frac{2}{m}\right)\right)$. Combining the diagonal and off-diagonal terms we get \begin{align*}
	\mathbb{E}\left[ \|I-2\eta A^*A\|_F^2\right] = d\left( 1- 4\eta + 4\eta^2 +\frac{8\eta^2}{m}\right) + (d^2-d)\frac{4\eta^2}{m}=4d\left[ \left(\frac{d+m+1}{m}\right)\eta^2-\eta+\frac{1}{4}\right].
	\end{align*}
\end{proof}

In much the same way we now compute the expected Frobenius norm of $2\eta A^*A$. 
\begin{lemma}\label{gradient}
	Let $A$ be an $m\times d$ random matrix with entries $A_{i,j}$ selected i.i.d. from $N(0,\frac{1}{m})$ and let $\eta >0$. Then the squared Frobenius norm of $2\eta A^*A$ has expected size\begin{align*}
	\mathbb{E}\left[\|2\eta A^*A\|_F^2\right]=4d\left[ \left(\frac{d+m+1}{m}\right) \eta^2\right].
	\end{align*}
\end{lemma}

\begin{proof}
	We will compute the expected squared Frobenius norm for $A^*A$ and pick up the $4\eta^2$ by linearity. The $d^2-d$ off-diagonal entries of $A^*A$ are precisely $\langle A_i,A_j\rangle$ for $i$ and $j$ distinct. Thus, by Lemma \ref{inner_product_lemma} these terms each contribute $\frac{1}{m}$. Moreover, the diagonal terms are $\|A_i\|_2^2$ which have expected squared size $1+\frac{2}{m}$. Thus we have \begin{align*}
	\mathbb{E}\left[ \|A^*A\|_F^2\right] = (d^2-d)\left(\frac{1}{m}\right) + d\left( 1+\frac{2}{m}\right) = d\left[\frac{d+m+1}{m}\right].
	\end{align*}
\end{proof}

Finally, in order to use Lemma \ref{Frobenius_Lemma}, we need to verify that the matrices $2\eta A^*A$ and $I-2\eta A^*A$ satisfy its conditions. This lemma does precisely that. 
\begin{lemma}\label{randomness_verification}
	Let $A$ be an $m\times d$ random matrix selected i.i.d. from $N(0,\frac{1}{m})$. Then the matrices $2\eta A^*A$ and $I-2\eta A^*A$ each satisfy \begin{enumerate}
		\item The columns of each matrix have the same expected squared length, and 
		\item For $V=2\eta A^*A$ or $V=I-2\eta A^*A$, two distinct entries $V_{i,j}$ and $V_{i,k}$ from the same row of $V$ have $\mathbb{E}\left[ V_{i,j}V_{i,k}\right] =0$. 
	\end{enumerate}
\end{lemma}

\begin{proof}
	It is straightforward to see that both $I-2\eta A^*A$ and $2\eta A^*A$ have columns with the same expected squared length by their definition. Let us verify condition $(2)$ directly for each matrix. 
	
	Let us begin with $V=2\eta A^*A$. Two distinct entries from row $i$ of $V$ are either
	\begin{itemize} 
		\item $\langle A_i, A_j\rangle$ and $\langle A_i, A_k\rangle$ for $j\neq k$ and neither equal to $i$, or 
		\item $\langle A_i, A_j\rangle$ and $\|A_i\|_2^2$ for $j\neq i$.
	\end{itemize}
In both of these situations, the expectation of their products are handled by equation \eqref{different_inner_products} of Lemma \ref{inner_product_lemma} taking first $p=k$ then $p=i$. 

Now, let $V=I-2\eta A^*A$. Akin to when $V=2\eta A^*A$, two distinct entries from a row of $V$ look like either \begin{itemize}
	\item $-2\eta \langle A_i, A_j\rangle$ and $-2\eta \langle A_i,A_k\rangle$ for $j\neq k$ and neither equal to $i$, or 
	\item $1-2\eta \|A_i\|_2^2$ and $-2\eta \langle A_i, A_j\rangle$  for $j\neq i$. 
\end{itemize}
Again, the first of these two possibilities is handled by equation \eqref{different_inner_products} of Lemma \ref{inner_product_lemma}. Now all that is left to observe is that: \begin{align*}
\mathbb{E}\left[ \left(1-2\eta \|A_i\|_2^2\right)\cdot  \langle A_i, A_j\rangle \right] = \mathbb{E}\left[ \langle A_i, A_j\rangle \right] - \mathbb{E}\left[2\eta \|A_i\|_2^2 \langle A_i, A_j\rangle\right]=0.
\end{align*}
The first term is $0$ by a simple expansion while the second term is zero again by Lemma \ref{inner_product_lemma}. 

\end{proof}

Finally we can combine all the previous results from this section to be able to say that in the average case, look ahead thresholding outperforms hard thresholding on arbitrary vectors.

\begin{theorem}\label{better_projection_theorem}
	Fix vectors $z,x^*$ in $\mathbb{R}^d$ with $x^*$ $s-$sparse and let $A$ be an $m\times d$ random matrix with i.i.d. $N(0,\frac{1}{m})$ entries. Then, \begin{align}
	\mathbb{E}\left[ \|H_{s,\eta} (z)-x^*\|_2\right] \leq \rho \|z-x^*\|_2 \label{better_thresholding_inequality}
	\end{align}
	for $\rho<2$ whenever $0 < \eta < \frac{1}{2} \frac{m}{m+d+1}$. 

\end{theorem}
\begin{proof}
	From Lemma \ref{transformation_lemma}, we have that \begin{align}
	\mathbb{E}\left[\|H_{s,\eta}(z)-x^*\|_2 \right]\leq\mathbb{E}\left[ \sqrt{\|z-\ell_\eta\|_2^2+\|\ell_\eta-x^*\|_2^2}+\|\ell_\eta - x^*\|_2\right]\label{eta_projection_expansion}.
	\end{align}
	First we will deal with the square root term. By substituting from Lemmas \ref{I_mins} and \ref{gradient} and using Jensen's inequality we get \begin{align}
	\mathbb{E}\left[\sqrt{\|z-\ell_\eta\|_2^2+\|\ell_\eta-x^*\|_2^2} \right]&\leq \left( \mathbb{E}\left[ \|z-\ell_\eta\|_2^2+\|\ell_\eta-x^*\|_2^2\right]\right)^{1/2}\nonumber\\
	&=\left( \mathbb{E} \left[ \|2\eta A^*A\|_F^2 + \|I-2\eta A^*A\|_F^2\right]\right)^{1/2}\cdot \|z-x^*\|_2\nonumber\\
	&= \left( \frac{8(d+m+1)}{m}\eta^2-4\eta +1\right)^{1/2}\cdot \|z-x^*\|_2\label{square_root_part}.
	\end{align}
	Now the second term of \eqref{eta_projection_expansion} can be handled with just Lemma \ref{gradient} to yield:  \begin{align}
	\mathbb{E}\left[ \|\ell_\eta - x^*\|_2 \right] &\leq \left(\mathbb{E}\left[\|(I-A^*A)(x^*-z)\|_2^2\right] \right)^{1/2}\nonumber\\
	&=\left(\frac{4(d+m+1)}{m}\eta^2-4\eta +1\right)^{1/2}\|z-x^*\|_2\label{second_term}.
	\end{align}
	Noticing that \eqref{second_term} is dominated by \eqref{square_root_part} for every value of $\eta$, we can make the greatly simplifying substitution that \begin{align*}
	\mathbb{E}\left[ \|H_{s,\eta} (z)-x^*\|_2\right]\leq 2\left(\frac{8(d+m+1)}{m}\eta^2 - 4\eta +1\right)^{1/2}\|z-x^*\|_2\label{simplified bound}
	\end{align*}
and it is easy to verify that $\frac{8(d+m+1)}{m}\eta^2 -4\eta +1$ is less than $1$ for precisely the specified values of eta.
\end{proof}

	We have two comments related to the analysis in this section. First, after equation \eqref{second_term} we make a simplifying substitution to reduce the complexity greatly. If desired, we could use numerical software to expand the range of $\eta$ values that give better average case performance.

	Second, the result in Theorem \ref{better_projection_theorem} is an attempt to say that, when compared to hard thresholding, look ahead thresholding returns sparse vectors that are closer to the desired solution on average. However, when $H_{s,\eta}$ is used in practice, e.g. as a step in ILAT, the points to threshold are far from arbitrary. In fact, they rely quite heavily on the draw of the random matrix $A$. Perhaps surprisingly, the experimental results of the following section seem to suggest that this interaction actually \textit{improves} the performance of look ahead thresholding relative to hard thresholding. That is, the range of values of $\eta$ for which ILAT outperforms IHT is much larger than predicted by Theorem \ref{better_projection_theorem}. Certainly there is more to be understood here.

\section{Experiments}\label{experiments_section}
While the last section was theoretical justification for the use of look ahead thresholding over hard thresholding, this section is experimental justification. 

 We perform three main experiments to support our claims. First, we show how well ILAT performs relative to IHT when the algorithms' iterations are held constant. Second, because ILAT takes extra computational power, we show that even when holding computations approximately constant we outperform IHT. Finally, we validate directly Theorem \ref{better_projection_theorem} by testing the thresholders' performances outside the framework of iterative algorithms. 

\begin{figure}
	\centering
	\includegraphics[width=0.7\linewidth]{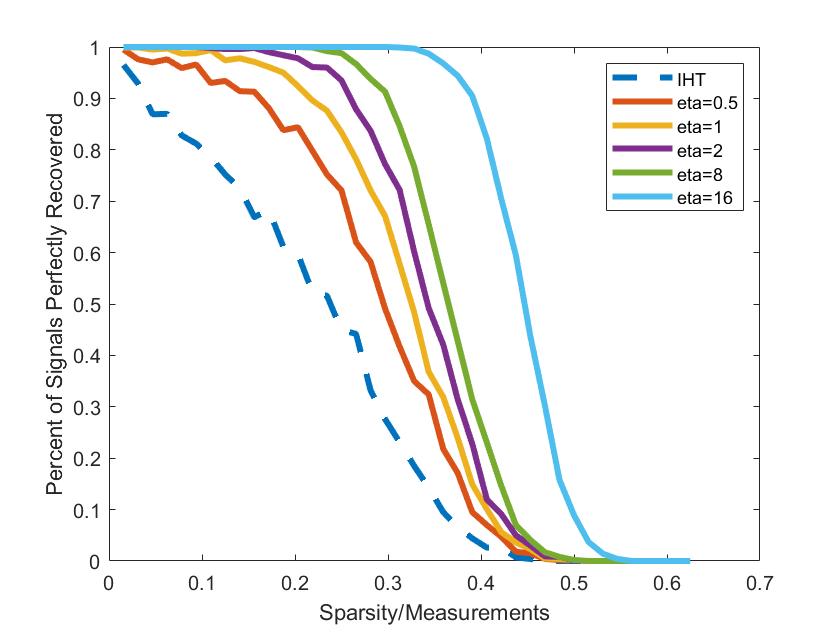}
	\caption{Percentage of signals perfectly reconstructed using IHT (dashed) and ILAT (solid, $\eta$ increases to the right) using the same number of iterations}
	\label{fig:betterplotnoiseless}
\end{figure}

For our first experiment, we record the percentage of exact recoveries for a range of sparsity values for both IHT and ILAT with a variety of $\eta$ values. The measurement matrix we use is $128\times 256$ i.i.d. Gaussian scaled to have operator norm $1$. Figure \ref{fig:betterplotnoiseless} shows ILAT outperforming IHT for a range of $\eta$ values.

\begin{figure}
	\centering
	\begin{subfigure}{.5\textwidth}
		\centering
		\includegraphics[width=\linewidth]{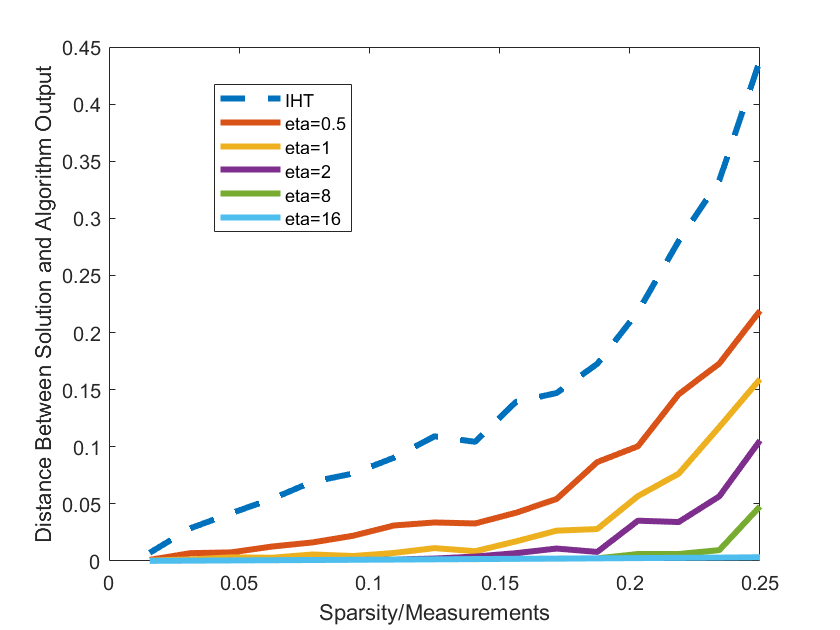}
		\caption{Error with 30dB SNR}
		\label{fig:sub1a}
	\end{subfigure}%
	\begin{subfigure}{.5\textwidth}
		\centering
		\includegraphics[width=\linewidth]{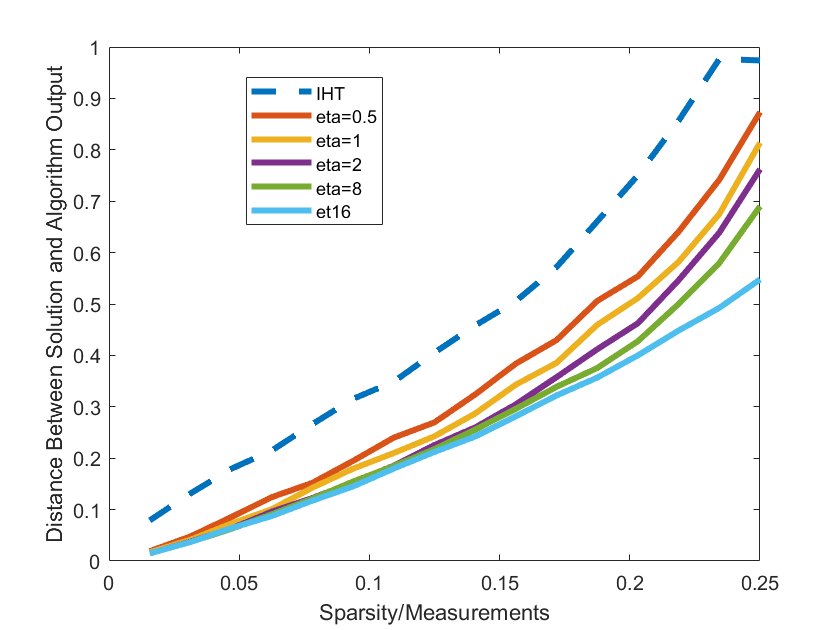}
		\caption{Error with 0dB SNR}
		\label{fig:sub1b}
	\end{subfigure}
	\caption{Reconstruction Error of IHT (dashed) and ILAT (solid, $\eta$ increasing to the right)}
	\label{fig:noisyfig1}
\end{figure}

We also tested our algorithm in the case where our measurements are corrupted by random noise. The performance metric used here is Euclidean distance between the true solution and the algorithm's output after $1000$ iterations. We performed these experiments in the same settings as the above example. For Figure \ref{fig:sub1a}, our measurements had a $30$ dB SNR and in Figure \ref{fig:sub1b} a 0 dB SNR. As above, Figure \ref{fig:noisyfig1} show ILAT achieves smaller error than IHT.

The two previous experiments held the number of iterations constant. Since look ahead thresholding is more costly than hard thresholding, we need to ask whether ILAT is superior to IHT when we hold constant the amount of computations or run-time. Notice that the expensive step is the gradient computation and that ILAT requires two gradients per iteration while IHT requires only one. With this in mind, we compared the percentage of perfect recoveries for IHT against ILAT when IHT is allowed twice the number of iterations.

\begin{figure}[H]
	\centering
	\begin{subfigure}{.5\textwidth}
		\centering
		\includegraphics[width=\linewidth]{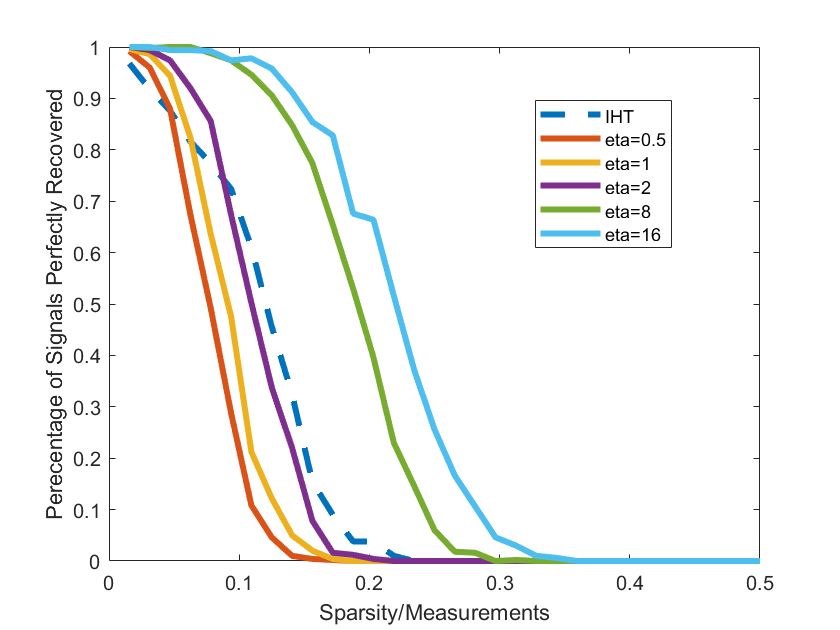}
		\caption{100 Gradient Computations}
		\label{fig:const_time1}
	\end{subfigure}%
	\begin{subfigure}{.5\textwidth}
		\centering
		\includegraphics[width=\linewidth]{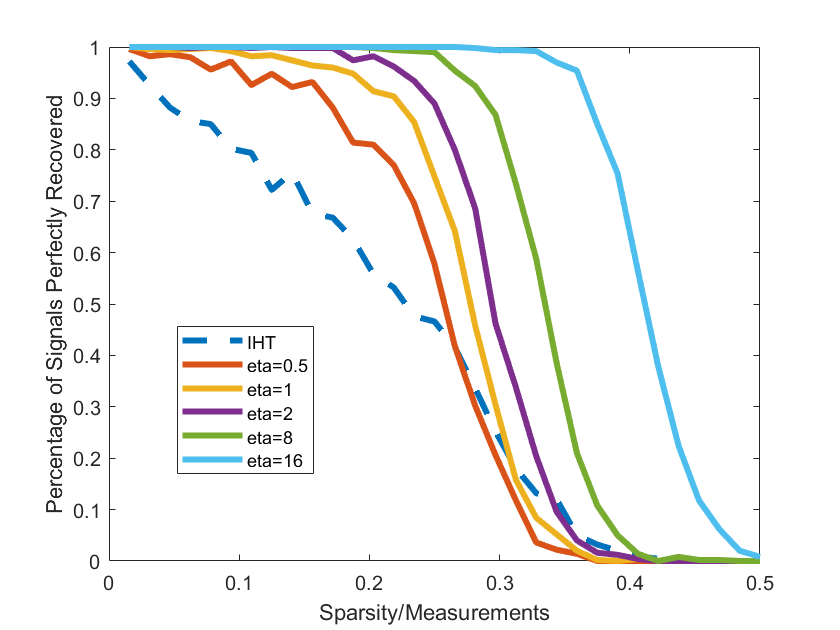}
		\caption{500 Gradient Computations}
		\label{fig:constant_time2}
	\end{subfigure}
	\caption{Percentage of signals perfectly reconstructed using IHT (dashed) and ILAT (solid, $\eta$ increasing to the right). Here, IHT is allowed twice as many iterations as ILAT}
	\label{fig:constant_time}
\end{figure}

On the left in Figure \ref{fig:const_time1}, IHT is given 100 iterations for recovery while ILAT is given 50 iterations for recovery. On the right in Figure \ref{fig:constant_time2}, the algorithms are given 500 and 250 iterations, respectively. Again, the entirety of Figure \ref{fig:constant_time} shows ILAT exceeding the recovery ability of IHT for appropriately chosen values of $\eta$. 

For our final experiment, we compare look ahead thresholding directly to hard thresholding when they are not simply a step in another algorithm. This lets us directly validate the conclusions of Theorem \ref{better_projection_theorem}. 

In the following, we generate a random sparse signal $x^*$ where the support is chosen uniformly at random and the entries on that support are i.i.d. $N(0,1)$. Morevoer, we pick a dense vector $z$ from $N(0,I)$ and the matrix $A$ has entries chosen i.i.d. from $N(0,1)$ then normalized to have $\|A\|_{op}=1$. Then, we measure how close the thresholdings $H_s(z)$ and $H_{s,\eta}(z)$ are to $x^*$. 

In light of the statement of \ref{better_projection_theorem}, Figure \ref{fig:arbitrary_vector} below shows the plots for $\|H_{s,\eta} (z)-x^*\|_2/\|z-x^*\|_2$ and $\|H_s(z)-x^*\|_2/\|z-x^*\|_2$. Note: here we only use $\eta=0.5,1$ because for $\eta=2,8,16$, the performance is indistinguishable from the $\eta=1$ case. Finally, on the left we use a $128\times 256$ matrix while the right figure uses only a $64\times 256$ measurement matrix to define $H_{s,\eta}$.

\begin{figure}[H]
	\centering
	\begin{subfigure}{.5\textwidth}
		\centering
		\includegraphics[width=\linewidth]{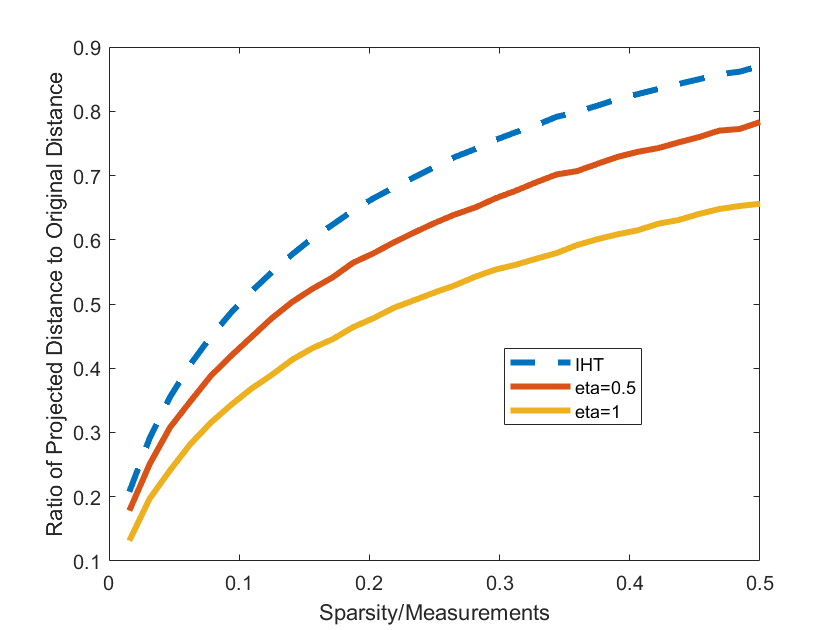}
		\caption{Thresholding Error with 128 Measurements}
		\label{fig:sub1}
	\end{subfigure}%
	\begin{subfigure}{.5\textwidth}
		\centering
		\includegraphics[width=\linewidth]{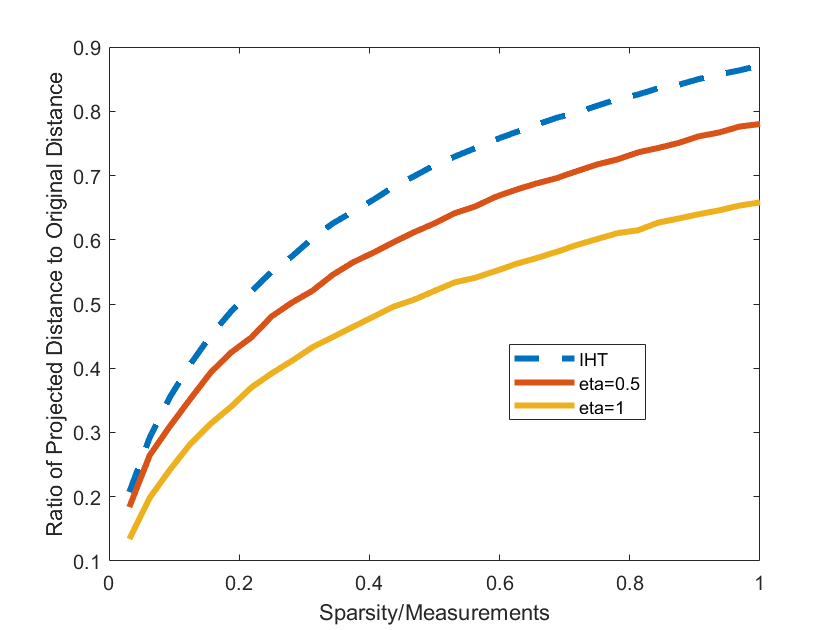}
		\caption{Thresholding Error with 64 Measurements}
		\label{fig:sub2}
	\end{subfigure}
	\caption{Distance from output $H_s(z)$ (dashed) and $H_{s,\eta} (z)$ (solid, $\eta$ increasing to the right) to solution $x^*$ divided by distance from $z$ to $x^*$}
	\label{fig:arbitrary_vector}
\end{figure}

\section{Conclusions}\label{conclusions_section}
The main contribution of this paper is the introduction of an alternative thresholding rule to hard thresholding for compressed sensing. Our thresholder uses specific problem instance information, namely the vector of measurements and measurement matrix, to help choose which entries to retain. This stands in stark contrast to hard thresholding which makes no use of such information and, as we have shown, is therefore suboptimal. 

Our thresholding rule is meant to act as a tool for use throughout compressed sensing. We investigated in-depth one of these cases by defining Iterative Look Ahead Thresholding, a variation on Iterative Hard Thresholding, which uses our new thresholding rule. We showed that the worst case performance of ILAT is comparable to IHT while the average case performance and experimental results exhibit greatly improved signal recovery. Even better, though our thresholding rule requires extra computational complexity, the amount of additional time is small and we have shown that even when computations are held approximately constant ILAT still excels.

There are still outstanding questions related to look ahead thresholding. First, our theoretical analysis of Section \ref{average_section} and experimental results of Section \ref{experiments_section} both suggest improved performance relative to hard thresholding. However, the experimental results are stronger than the theoretical results, i.e. look ahead thresholding works for a much larger range of $\eta$ values than predicted. Moreover, after explaining this disparity, we would like to understand how to pick the optimal value of $\eta$ given the ambient dimension, number of measurements, sparsity level, and perhaps a signal model.

In addition to the questions we have about ILAT, we also encourage others to explore the utility of look ahead thresholding in a variety of different CS algorithms. In particular, we are aware that CoSaMP makes critical usage of a thresholding step so perhaps a more accurate thresholding rule could increase the algorithm's accuracy.

\section*{Acknowledgements}
The Author would like to thank Alexander Powell and Zack Tripp for insights, advice, and ideas without which this paper would have been impossible.

\bibliographystyle{amsplain}
\bibliography{biblio}

\end{document}